\documentclass[oneside,reqno,english]{amsart}
\usepackage[T1]{fontenc}
\usepackage[utf8]{inputenc}
\usepackage{xcolor}
\usepackage{babel}
\usepackage{prettyref}
\usepackage{mathtools}
\usepackage{amstext}
\usepackage{amsthm}
\usepackage{amssymb}
\usepackage[pdfusetitle,
 bookmarks=true,bookmarksnumbered=false,bookmarksopen=false,
 breaklinks=false,pdfborder={0 0 0},pdfborderstyle={},backref=false,colorlinks=false]
 {hyperref}
\hypersetup{
 colorlinks=true,citecolor=blue,linkcolor=blue,linktocpage=true}

\makeatletter
\numberwithin{equation}{section}
\numberwithin{figure}{section}

\usepackage{prettyref}

\newrefformat{cor}{Corollary~\ref{#1}}
\newrefformat{subsec}{Section~\ref{#1}}
\newrefformat{lem}{Lemma~\ref{#1}}
\newrefformat{thm}{Theorem~\ref{#1}}
\newrefformat{sec}{Section~\ref{#1}}
\newrefformat{chap}{Chapter~\ref{#1}}
\newrefformat{prop}{Proposition~\ref{#1}}
\newrefformat{exa}{Example~\ref{#1}}
\newrefformat{tab}{Table~\ref{#1}}
\newrefformat{rem}{Remark~\ref{#1}}
\newrefformat{def}{Definition~\ref{#1}}
\newrefformat{fig}{Figure~\ref{#1}}
\newrefformat{claim}{Claim~\ref{#1}}
\newrefformat{assu}{Assumption~\ref{#1}}

\makeatother

\theoremstyle{plain}
\newtheorem{thm}{\protect\theoremname}[section]
\newtheorem{prop}[thm]{\protect\propositionname}
\theoremstyle{definition}
\newtheorem{defn}[thm]{\protect\definitionname}
\theoremstyle{plain}
\newtheorem{lem}[thm]{\protect\lemmaname}
\newtheorem{cor}[thm]{\protect\corollaryname}
\theoremstyle{remark}
\newtheorem{rem}[thm]{\protect\remarkname}
\theoremstyle{definition}
\newtheorem{example}[thm]{\protect\examplename}
\providecommand{\corollaryname}{Corollary}
\providecommand{\definitionname}{Definition}
\providecommand{\examplename}{Example}
\providecommand{\lemmaname}{Lemma}
\providecommand{\propositionname}{Proposition}
\providecommand{\remarkname}{Remark}
\providecommand{\theoremname}{Theorem}

\begin{document}
\title{Kernel Radon-Nikodym Derivatives for Random Matrix Products}
\begin{abstract}
This paper studies kernel Radon-Nikodym derivatives for the one-step
shift of time-indexed positive definite kernels associated with random
matrix products. The problem is to determine when the shifted kernel
is dominated by the original kernel and to identify the corresponding
Radon-Nikodym derivative. We treat two concrete classes of multiplicative
walks: ensembles with inhomogeneous variances and Gaussian Kraus products.
In both settings, the shifted kernel inequality reduces to a one-step
condition on the diagonal moments, and the Radon-Nikodym derivative
is described explicitly by a fiberwise sequence in the time variable.
In the inhomogeneous variance model, the diagonal compression is governed
by a nonnegative matrix $S$, which yields an explicit coordinate
formula for the fibers. In the Gaussian Kraus model, the diagonal
moments are generated by a completely positive map $\Psi$, and the
shifted kernel inequality is equivalent to the condition ${\Psi\left(I\right)\le I}$.
\end{abstract}

\author{James Tian}
\address{Mathematical Reviews, 535 W. William St, Suite 210, Ann Arbor, MI
48103, USA}
\email{james.ftian@gmail.com}

\maketitle
\tableofcontents{}

\section{Introduction}\label{sec:1}

We work with operator-valued positive definite kernels on the free
semigroup $\mathbb{F}^{+}_{d}$. Writing $\Sigma$ for the right shift,
a kernel $K$ gives a shifted kernel $K_{\Sigma}$, and we consider
the order inequality $K_{\Sigma}\le K$. When this holds, the Kolmogorov
decomposition of $K$ identifies a canonical contractive operator
on the dilation space implementing the shift, and this operator is
identified as the Radon-Nikodym derivative $dK_{\Sigma}/dK$ in the
cone of positive definite kernels.

In the time-diagonal situations arising from the random matrix models
below, this Radon-Nikodym derivative breaks into a concrete sequence
of one-step fibers $D_{m}$, and these fibers measure the change from
depth $m$ to depth $m+1$ in an intrinsic way. In particular, the
fibers are attached to the kernel inequality itself, rather than chosen
afterward as normalizations of successive moments.

The kernel background used in \prettyref{sec:2} is standard: reproducing
kernel and operator-valued positive definite kernel methods \cite{MR51437,MR4439542,MR4914974},
together with the free semigroup shift and its role in dilation theory
for row contractions \cite{MR671311,MR744917,MR769375,MR1407612,MR1618326,MR1681749,MR2225440,Tian:2026aa}.
There is also a large literature on Radon-Nikodym type results for
completely positive maps and quantum operations \cite{MR932932,MR1608473,MR2014842,MR2761319,MR3697039,MR247483,MR394232,MR4574232,MR5002158}.
The point here is different. We do not compare two completely positive
maps directly. Instead, we compare a kernel with its shift and take
the Radon-Nikodym derivative inside the cone of positive definite
kernels. 

We study this in two classes of multiplicative walks. The model with
inhomogeneous variances considered below also sits against the broader
background of random matrix theory and products of random matrices
\cite{MR2760897,MR2567175,MR2653641,MR3901642,MR4552960,MR4664588,MR4870344}.
In particular, the diagonal compression is governed by a nonnegative
matrix $S$, and the fiberwise densities admit an explicit coordinate
description in terms of the vectors $S^{m}\mathbf{1}$. This gives
a finite-dimensional setting in which the kernel Radon-Nikodym fibers
can be written down completely. 

In the Gaussian Kraus product setting, the diagonal moments are governed
by a completely positive map $\Psi$, placing the model alongside
the operator-algebraic literature on completely positive maps, extension,
and kernel domination \cite{MR394232,MR4574232,MR5002158}, and the
shift domination problem becomes an exact one-step criterion:
\[
K_{\Sigma}\le K\iff\Psi\left(I\right)\le I.
\]
At the same time, the kernel Radon-Nikodym construction identifies
the corresponding one-step fibers $D_{m}$ explicitly. This gives
a concrete fiberwise description of the one-step comparison in the
model and shows that the kernel order can be read directly from the
one-step action of $\Psi$. We also give an example showing that the
condition $\Psi\left(I\right)\le I$ is independent from the spectral
radius condition $r\left(\Psi\right)<1$.

In this paper, we have restricted attention to the one-parameter shift
relevant to the studied models. Extensions to genuinely branching
$\mathbb{F}^{+}_{d}$ ($d>1$) settings, with corresponding multishift
Radon-Nikodym comparisons, appear to require additional structure
and are left for future work.

The paper is organized as follows. In \prettyref{sec:2} we collect
the kernel preliminaries, including the Kolmogorov decomposition and
the Radon-Nikodym derivative for positive definite kernels, and we
introduce the shifted kernel on $\mathbb{F}^{+}_{d}$. In \prettyref{sec:3}
we study diagonal compressions coming from ensembles with inhomogeneous
variances and write down the resulting fiberwise densities. In \prettyref{sec:4}
we study the Gaussian Kraus product model, prove the criterion $\Psi\left(I\right)\le I$,
and give an example separating $\Psi\left(I\right)\le I$ from the
spectral radius condition $r\left(\Psi\right)<1$.

\section{Preliminaries}\label{sec:2}

In this section we collect the basic definitions and facts used later.
In particular, we recall the Kolmogorov decomposition and the Radon-Nikodym
derivative for positive definite kernels, and we introduce the shifted
kernel on the free semigroup $\mathbb{F}^{+}_{d}$. 

For background on positive definite and operator-valued positive definite
kernels, we refer to \cite{MR51437,MR1454031,MR2265340,MR2603770,MR2938971,MR4439542}. 

Let $H$ be a Hilbert space, and $\mathcal{L}\left(H\right)$ the
space of bounded linear operators. Throughout, all inner products
are linear in the second variable. 

A kernel $K:X\times X\to\mathcal{L}\left(H\right)$ is called positive
definite if for every finite family $\left\{ \left(x_{j},u_{j}\right)\right\} ^{n}_{j=1}\subset X\times H$,
\[
\sum^{n}_{i,j=1}\left\langle u_{i},K\left(x_{i},x_{j}\right)u_{j}\right\rangle \ge0.
\]
By the Kolmogorov decomposition theorem, every positive definite kernel
$K$ admits a Hilbert space $\mathcal{K}_{K}$ and a map $V_{K}:X\to\mathcal{L}\left(H,\mathcal{K}_{K}\right)$
such that 
\[
K\left(x,y\right)=V_{K}\left(x\right)^{*}V_{K}\left(y\right),\qquad x,y\in X,
\]
and 
\[
\overline{span}\left\{ V_{K}\left(x\right)h:x\in X,\ h\in H\right\} =\mathcal{K}_{K}.
\]
We refer to $\left(V_{K},\mathcal{K}_{K}\right)$ as a minimal Kolmogorov
decomposition of $K$. 

Fix $d\in\mathbb{N}$. Let $\mathbb{F}{}^{+}_{d}$ be the free semigroup
on generators $\left\{ 1,\dots,d\right\} $, with neutral element
$\emptyset$. Let $B\subseteq\mathcal{L}\left(H\right)$ be a unital
$C^{*}$-subalgebra, and let 
\[
K:\mathbb{F}^{+}_{d}\times\mathbb{F}^{+}_{d}\to B
\]
be a $B$-valued positive definite kernel. We define the shifted kernel
by 
\[
K_{\Sigma}\left(\alpha,\beta\right)=\sum^{d}_{i=1}K\left(\alpha i,\beta i\right).
\]

\begin{prop}
\label{prop:1-3} Let $X$ be a set, and let $K,L:X\times X\to\mathcal{L}\left(H\right)$
be positive definite kernels with $K\le L$. Let 
\[
L\left(x,y\right)=V_{L}\left(x\right)^{*}V_{L}\left(y\right)
\]
be a minimal Kolmogorov decomposition of $L$ on $\mathcal{K}_{L}$.
Then there exists a unique operator $D\in\mathcal{L}\left(\mathcal{K}_{L}\right)$
such that 
\begin{equation}
0\le D\le I,\qquad K\left(x,y\right)=V_{L}\left(x\right)^{*}DV_{L}\left(y\right),\qquad x,y\in X.\label{eq:2-1}
\end{equation}
\end{prop}

\begin{proof}
This is standard in the CP setting (see e.g., \cite{MR247483,MR1608473}),
but for completeness we include a brief sketch of the proof. 

Choose a minimal Kolmogorov decomposition $K\left(x,y\right)=V_{K}\left(x\right)^{*}V_{K}\left(y\right)$
of $K$ on a Hilbert space $\mathcal{K}_{K}$. Set 
\[
\mathcal{D}_{L}=span\left\{ V_{L}\left(x\right)h:x\in X,\ h\in H\right\} \subseteq\mathcal{K}_{L}.
\]
Define $T_{0}:\mathcal{D}_{L}\to\mathcal{K}_{K}$ on generators by
\[
T_{0}\left(\sum\nolimits^{n}_{j=1}V_{L}\left(x_{j}\right)h_{j}\right)=\sum\nolimits^{n}_{j=1}V_{K}\left(x_{j}\right)h_{j}.
\]
We first show that $T_{0}$ is well defined. Suppose $\sum^{n}_{j=1}V_{L}\left(x_{j}\right)h_{j}=0$.
Then 
\begin{align*}
\left\Vert \sum\nolimits^{n}_{j=1}V_{K}\left(x_{j}\right)h_{j}\right\Vert ^{2} & =\sum\nolimits^{n}_{i,j=1}\left\langle h_{i},K\left(x_{i},x_{j}\right)h_{j}\right\rangle \\
 & \le\sum^{n}_{i,j=1}\left\langle h_{i},L\left(x_{i},x_{j}\right)h_{j}\right\rangle \\
 & =\left\Vert \sum\nolimits^{n}_{j=1}V_{L}\left(x_{j}\right)h_{j}\right\Vert ^{2}=0,
\end{align*}
so $\sum^{n}_{j=1}V_{K}\left(x_{j}\right)h_{j}=0$. Thus $T_{0}$
is well defined.

The same computation shows that $T_{0}$ is contractive on $\mathcal{D}_{L}$:
for $\xi=\sum^{n}_{j=1}V_{L}\left(x_{j}\right)h_{j}$ we have 
\[
\left\Vert T_{0}\xi\right\Vert ^{2}=\sum^{n}_{i,j=1}\left\langle h_{i},K\left(x_{i},x_{j}\right)h_{j}\right\rangle \le\sum^{n}_{i,j=1}\left\langle h_{i},L\left(x_{i},x_{j}\right)h_{j}\right\rangle =\left\Vert \xi\right\Vert ^{2}.
\]
Hence $T_{0}$ extends uniquely to a contraction $W:\mathcal{K}_{L}\to\mathcal{K}_{K}$.
By construction, $WV_{L}\left(x\right)=V_{K}\left(x\right)$, $x\in X$.
Set $D=W^{*}W$. Then $0\le D\le I$. For $x,y\in X$ and $u,v\in H$,
\begin{align*}
\left\langle u,K\left(x,y\right)v\right\rangle  & =\left\langle V_{K}\left(x\right)u,V_{K}\left(y\right)v\right\rangle \\
 & =\left\langle WV_{L}\left(x\right)u,WV_{L}\left(y\right)v\right\rangle \\
 & =\left\langle V_{L}\left(x\right)u,W^{*}WV_{L}\left(y\right)v\right\rangle \\
 & =\left\langle u,V_{L}\left(x\right)^{*}DV_{L}\left(y\right)v\right\rangle .
\end{align*}
Therefore $K\left(x,y\right)=V_{L}\left(x\right)^{*}DV_{L}\left(y\right)$. 

For uniqueness, suppose $D_{1},D_{2}\in\mathcal{L}\left(\mathcal{K}_{L}\right)$
satisfy 
\[
V_{L}\left(x\right)^{*}D_{1}V_{L}\left(y\right)=V_{L}\left(x\right)^{*}D_{2}V_{L}\left(y\right)\qquad x,y\in X.
\]
Then 
\[
\left\langle V_{L}\left(x\right)u,\left(D_{1}-D_{2}\right)V_{L}\left(y\right)v\right\rangle =0
\]
for all $x,y,u,v$. Since $span\left\{ V_{L}\left(x\right)h:x\in X,\ h\in H\right\} $
is dense in $\mathcal{K}_{L}$, it follows that $D_{1}=D_{2}$.
\end{proof}

\begin{defn}
We denote the operator in \prettyref{eq:2-1} by 
\[
\frac{dK}{dL}=D
\]
and refer to it as the Radon-Nikodym derivative of $K$ with respect
to $L$.
\end{defn}

We now specialize to the shifted kernel on $\mathbb{F}^{+}_{d}$.
Let $K:\mathbb{F}^{+}_{d}\times\mathbb{F}^{+}_{d}\to\mathcal{L}\left(H\right)$
be positive definite, with minimal Kolmogorov decomposition 
\[
K\left(\alpha,\beta\right)=V\left(\alpha\right)^{*}V\left(\beta\right)
\]
on a Hilbert space $\mathcal{K}$. Set 
\[
\mathcal{D}=span\left\{ V\left(\alpha\right)h:\alpha\in\mathbb{F}^{+}_{d},\ h\in H\right\} .
\]

\begin{lem}
\label{lem:1-4} The following are equivalent: 
\begin{enumerate}
\item $K_{\Sigma}\le K$; 
\item there exists a contractive operator $T:\mathcal{K}\to\mathcal{K}^{\oplus d}$
such that 
\[
TV\left(\alpha\right)h=\left(V\left(\alpha1\right)h,\dots,V\left(\alpha d\right)h\right),\qquad\alpha\in\mathbb{F}^{+}_{d},\ h\in H.
\]
\end{enumerate}
When these conditions hold, if $P_{i}:\mathcal{K}^{\oplus d}\to\mathcal{K}$
denotes the $i$-th coordinate projection and $T_{i}=P_{i}T$, then
\[
T_{i}V\left(\alpha\right)h=V\left(\alpha i\right)h,\qquad\sum^{d}_{i=1}T^{*}_{i}T_{i}\le I_{\mathcal{K}}.
\]
Moreover, the operator $T$, and hence each $T_{i}$, is uniquely
determined. 
\end{lem}

\begin{proof}
Assume first that $K_{\Sigma}\le K$. Define $T_{0}:\mathcal{D}\to\mathcal{K}^{\oplus d}$
by 
\[
T_{0}\left(\sum\nolimits_{\alpha}V\left(\alpha\right)u_{\alpha}\right)=\left(\sum\nolimits_{\alpha}V\left(\alpha1\right)u_{\alpha},\dots,\sum\nolimits_{\alpha}V\left(\alpha d\right)u_{\alpha}\right),
\]
for finite sums. We show that $T_{0}$ is well defined. If $x=\sum\nolimits_{\alpha}V\left(\alpha\right)u_{\alpha}=0$,
then 
\begin{align*}
\left\Vert T_{0}x\right\Vert ^{2} & =\sum\nolimits^{d}_{i=1}\left\Vert \sum\nolimits_{\alpha}V\left(\alpha i\right)u_{\alpha}\right\Vert ^{2}\\
 & =\sum\nolimits_{\alpha,\beta}\left\langle u_{\alpha},\left(\sum\nolimits^{d}_{i=1}K\left(\alpha i,\beta i\right)\right)u_{\beta}\right\rangle \\
 & \le\sum\nolimits_{\alpha,\beta}\left\langle u_{\alpha},K\left(\alpha,\beta\right)u_{\beta}\right\rangle =\left\Vert x\right\Vert ^{2}=0.
\end{align*}
Hence $T_{0}x=0$, so $T_{0}$ is well defined. The same computation
shows that $T_{0}$ is contractive on $\mathcal{D}$. Since $\mathcal{D}$
is dense in $\mathcal{K}$, $T_{0}$ extends uniquely to a contractive
operator $T:\mathcal{K}\to\mathcal{K}^{\oplus d}$. By construction,
$TV\left(\alpha\right)h=\left(V\left(\alpha1\right)h,\dots,V\left(\alpha d\right)h\right)$. 

Conversely, assume such a contractive operator $T$ exists. Write
$T_{i}=P_{i}T$. Then for any finitely supported family $\left\{ u_{\alpha}\right\} \subset H$,
if $x=\sum_{\alpha}V\left(\alpha\right)u_{\alpha}$, we have 
\begin{align*}
\sum\nolimits_{\alpha,\beta}\left\langle u_{\alpha},\left(K\left(\alpha,\beta\right)-\sum\nolimits^{d}_{i=1}K\left(\alpha i,\beta i\right)\right)u_{\beta}\right\rangle  & =\left\Vert x\right\Vert ^{2}-\sum^{d}_{i=1}\left\Vert T_{i}x\right\Vert ^{2}\\
 & =\left\Vert x\right\Vert ^{2}-\left\Vert Tx\right\Vert ^{2}\ge0.
\end{align*}
Thus $K-K_{\Sigma}$ is positive definite, i.e. $K_{\Sigma}\le K.$

Finally, since $T$ is contractive, $\sum^{d}_{i=1}T^{*}_{i}T_{i}=T^{*}T\le I_{\mathcal{K}}$.
Uniqueness follows from the density of $\mathcal{D}$. 
\end{proof}

\begin{prop}
\label{prop:1-5} Under the hypotheses of \prettyref{lem:1-4}, 
\[
\frac{dK_{\Sigma}}{dK}=\sum^{d}_{i=1}T^{*}_{i}T_{i}.
\]
In particular, 
\[
K_{\Sigma}\le K\iff\sum^{d}_{i=1}T^{*}_{i}T_{i}\le I_{\mathcal{K}}.
\]
\end{prop}

\begin{proof}
For $\alpha,\beta\in\mathbb{F}^{+}_{d}$ and $u,v\in H$, 
\begin{align*}
\left\langle u,K_{\Sigma}\left(\alpha,\beta\right)v\right\rangle  & =\sum\nolimits^{d}_{i=1}\left\langle u,K\left(\alpha i,\beta i\right)v\right\rangle \\
 & =\sum\nolimits^{d}_{i=1}\left\langle V\left(\alpha i\right)u,V\left(\beta i\right)v\right\rangle \\
 & =\sum\nolimits^{d}_{i=1}\left\langle T_{i}V\left(\alpha\right)u,T_{i}V\left(\beta\right)v\right\rangle \\
 & =\left\langle V\left(\alpha\right)u,\left(\sum\nolimits^{d}_{i=1}T^{*}_{i}T_{i}\right)V\left(\beta\right)v\right\rangle .
\end{align*}
Therefore 
\[
K_{\Sigma}\left(\alpha,\beta\right)=V\left(\alpha\right)^{*}\left(\sum\nolimits^{d}_{i=1}T^{*}_{i}T_{i}\right)V\left(\beta\right).
\]
By \prettyref{prop:1-3}, the operator $\frac{dK_{\Sigma}}{dK}$ is
the unique positive operator on $\mathcal{K}$ implementing this factorization.
Hence 
\[
\frac{dK_{\Sigma}}{dK}=\sum^{d}_{i=1}T^{*}_{i}T_{i}.
\]
The final equivalence is exactly \prettyref{lem:1-4}. 
\end{proof}

We now pass from the abstract kernel setting to the random setting
that motivates the paper. In what follows, the kernels $K$ will be
built from moment expectations of multiplicative random walks, typically
of the form 
\[
K\left(m,n\right)=E_{\mathcal{B}}\mathbb{E}\left[X_{m}X^{*}_{n}\right],
\]
where $\left(X_{m}\right)$ is a product process and $E_{\mathcal{B}}$
is a conditional expectation onto a fixed $*$-subalgebra $\mathcal{B}$.
In these models the index set reduces to the time parameter $m\in\mathbb{N}$
(equivalently, the $d=1$ semigroup inside $\mathbb{F}^{+}_{d}$),
and the shift $\Sigma$ becomes the one-step time shift $m\mapsto m+1$.
The inequality $K_{\Sigma}\le K$ is then a concrete comparison between
successive moment levels, and the Radon-Nikodym derivative $\frac{dK_{\Sigma}}{dK}$
admits an explicit fiberwise description in terms of the one-step
evolution of the diagonal moments. 

The case $d>1$ is interesting in its own right, but it is not part
of the present scope.

\section{Gaussian products with inhomogeneous variances}\label{sec:3}

In this section we compute the Radon-Nikodym derivatives associated
with the time-shifted moment kernels of a multiplicative Gaussian
walk with inhomogeneous entry variances. After compressing to the
diagonal algebra, the kernel becomes diagonal in the time index and
satisfies an exact one-step recursion governed by $\Phi$, equivalently
by the linear operator $S$ on $\mathbb{R}^{N}$. This yields an explicit
criterion for the shifted kernel inequality and a fiberwise Radon-Nikodym
formula. We then study compressions onto block-diagonal subalgebras;
under a block-homogeneity assumption on $\sigma^{2}_{ij}$, one obtains
analogous criteria and formulas. 

Let $\left(A_{m}\right)_{m\ge1}$ be i.i.d. $N\times N$ random matrices
of the form 
\[
\left(A_{m}\right)_{ij}=\frac{\sigma_{ij}}{\sqrt{N}}g^{(m)}_{ij},
\]
where for each $m$ the variables $(g^{(m)}_{ij})_{i,j}$ are i.i.d.
complex Gaussians $\mathcal{N}_{\mathbb{C}}\left(0,1\right)$, and
the arrays for different $m$ are independent. Set 
\[
X_{0}\coloneqq I,\qquad X_{m}\coloneqq A_{m}A_{m-1}\cdots A_{1}\quad\text{for }m\ge1.
\]
We define the moment kernel 
\begin{equation}
K\left(m,n\right)\coloneqq\mathbb{E}\left[X_{m}\left(X_{n}\right)^{*}\right],\qquad m,n\in\mathbb{N},\label{eq:3-1}
\end{equation}
and, as before, 
\[
K_{\Sigma}\left(m,n\right)\coloneqq K\left(m+1,n+1\right).
\]

Let $D\subseteq M_{N}$ denote the diagonal algebra, and let $E_{D}:M_{N}\to D$
be the diagonal conditional expectation. Define the linear map $S:\mathbb{R}^{N}\to\mathbb{R}^{N}$
by 
\begin{equation}
\left(Sx\right)_{i}\coloneqq\sum^{N}_{j=1}\frac{\sigma^{2}_{ij}}{N}x_{j},\qquad x\in\mathbb{R}^{N},\label{eq:3-2}
\end{equation}
and define $\Phi:D\to D$ by 
\begin{equation}
\Phi\left(\mathrm{diag}\left(x\right)\right)\coloneqq\mathrm{diag}\left(Sx\right),\qquad x\in\mathbb{R}^{N}.\label{eq:3-3}
\end{equation}

We first isolate the diagonal part of the moment kernel and identify
its one-step evolution under multiplication by an independent increment.
\begin{lem}
\label{lem:3-1} Let the setting be as above, and $K$ as in \prettyref{eq:3-1}.
For all $m,n\in\mathbb{N}$ with $m\neq n$, one has $K\left(m,n\right)=0$.
If we set 
\[
M_{m}\coloneqq E_{D}K\left(m,m\right)=E_{D}\mathbb{E}\left[X_{m}X^{*}_{m}\right]\in D,
\]
then 
\begin{equation}
M_{m+1}=\Phi\left(M_{m}\right),\qquad M_{0}=I.\label{eq:3-4}
\end{equation}
Consequently, 
\begin{equation}
M_{m}=\Phi^{m}\left(I\right),\qquad\forall m\in\mathbb{N}.\label{eq:3-5}
\end{equation}
\end{lem}

\begin{proof}
Since each $A_{r}$ satisfies $A_{r}\stackrel{d}{=}e^{i\theta}A_{r}$,
we have $X_{m}\stackrel{d}{=}e^{im\theta}X_{m}$. Hence 
\[
K\left(m,n\right)=\mathbb{E}\left[X_{m}\left(X_{n}\right)^{*}\right]=e^{i\left(m-n\right)\theta}K\left(m,n\right)
\]
for every $\theta$, and therefore $K\left(m,n\right)=0$ when $m\neq n$.

For the recursion, write $X_{m+1}=A_{m+1}X_{m}$. By independence
of $A_{m+1}$ and $X_{m}$, 
\[
\mathbb{E}\left[X_{m+1}X^{*}_{m+1}\right]=\mathbb{E}\left[A_{m+1}X_{m}X^{*}_{m}A^{*}_{m+1}\right]=\mathbb{E}\left[\mathbb{E}\left[A_{m+1}X_{m}X^{*}_{m}A^{*}_{m+1}\mid X_{m}\right]\right].
\]
For a fixed matrix $Y$, a direct computation using the Gaussian second
moments gives 
\begin{equation}
E_{D}\mathbb{E}\left[A_{m+1}YA^{*}_{m+1}\right]=\Phi\left(E_{D}Y\right).\label{eq:3-6}
\end{equation}
Indeed, writing $A=A_{m+1}$ and taking diagonal entries, 
\[
\left(E_{D}\mathbb{E}\left[AYA^{*}\right]\right)_{ii}=\mathbb{E}\left[\sum\nolimits^{N}_{j,k=1}a_{ij}Y_{jk}\overline{a_{ik}}\right].
\]
By the Gaussian second moments (and independence across columns),
\[
\mathbb{E}\left[a_{ij}\overline{a_{ik}}\right]=0\ \ \text{for }j\ne k,\qquad\mathbb{E}[\left|a_{ij}\right|^{2}]=\frac{\sigma^{2}_{ij}}{N},
\]
so the sum collapses to 
\[
\left(E_{D}\mathbb{E}\left[AYA^{*}\right]\right)_{ii}=\sum^{N}_{j=1}\frac{\sigma^{2}_{ij}}{N}Y_{jj}=\left(\Phi\left(E_{D}Y\right)\right)_{ii}.
\]
Since both sides lie in $D$, this proves $E_{D}\mathbb{E}\left[AYA^{*}\right]=\Phi\left(E_{D}Y\right)$,
which is \prettyref{eq:3-6}. 

Applying this with $Y=X_{m}X^{*}_{m}$ and then taking expectation
yields 
\[
M_{m+1}=E_{D}\mathbb{E}\left[X_{m+1}X^{*}_{m+1}\right]=\Phi\left(E_{D}\mathbb{E}\left[X_{m}X^{*}_{m}\right]\right)=\Phi\left(M_{m}\right).
\]
Since $X_{0}=I$, we have $M_{0}=I$. Iterating \eqref{eq:3-4} gives
\eqref{eq:3-5}. 
\end{proof}

We next identify the exact condition under which the shifted kernel
inequality holds after diagonal compression.
\begin{thm}
\label{thm:3-2} Let $K_{D}\coloneqq E_{D}\circ K$. Write $r_{i}\coloneqq\sum^{N}_{j=1}\frac{\sigma^{2}_{ij}}{N}=\left(S\mathbf{1}\right)_{i}$,
$\mathbf{1}\coloneqq\left(1,\dots,1\right)^{\mathsf{T}}$, so that
\begin{equation}
M_{1}=E_{D}\mathbb{E}\left[A_{1}A^{*}_{1}\right]=\mathrm{diag}\left(r\right).\label{eq:3-7-1}
\end{equation}
Then the following are equivalent: 
\begin{enumerate}
\item $\left(K_{D}\right)_{\Sigma}\le K_{D}$. 
\item $M_{1}\le I$. 
\item $S\mathbf{1}\le\mathbf{1}$ entrywise, equivalently $r_{i}\le1$ for
all $i$. 
\end{enumerate}
Under these conditions, for all $m\in\mathbb{N}$, 
\[
K_{D}\left(m,m\right)=M_{m}=\Phi^{m}\left(I\right)=\mathrm{diag}\left(S^{m}\mathbf{1}\right),
\]
and 
\[
M_{m+1}\le M_{m}.
\]
\end{thm}

\begin{proof}
By \prettyref{lem:3-1}, the kernel $K_{D}$ is diagonal in $\left(m,n\right)$,
so the inequality $\left(K_{D}\right)_{\Sigma}\le K_{D}$ is equivalent
to 
\[
K_{D}\left(m+1,m+1\right)\le K_{D}\left(m,m\right)\qquad\text{for all }m,
\]
that is, 
\[
M_{m+1}\le M_{m}\qquad\text{for all }m.
\]
Since $\Phi$ is positive and order-preserving on $D$, the relation
$M_{m+1}=\Phi\left(M_{m}\right)$ from \prettyref{lem:3-1} gives
\[
M_{1}\le M_{0}=I\Longrightarrow M_{2}=\Phi\left(M_{1}\right)\le\Phi\left(I\right)=M_{1},
\]
and inductively 
\[
M_{m+1}\le M_{m}\qquad\text{for all }m.
\]
Thus (1) and (2) are equivalent.

From \prettyref{eq:3-7-1}, $M_{1}\le I$ is equivalent to $r\le\mathbf{1}$
entrywise, which is exactly $S\mathbf{1}\le\mathbf{1}$. This proves
the equivalence of (2) and (3). The formula for $K_{D}\left(m,m\right)$
follows from \eqref{eq:3-5}. 
\end{proof}

We now compute the Radon-Nikodym derivative for time-diagonal kernels
and then specialize the formula to $K_{D}$.
\begin{prop}
\label{prop:3-3} Let $H$ be a Hilbert space and let $K,L:\mathbb{N}\times\mathbb{N}\to\mathcal{L}\left(H\right)$
be positive definite kernels with $K\le L$. Assume that both kernels
are diagonal in the time index, i.e. 
\[
K\left(m,n\right)=0=L\left(m,n\right)\qquad\text{whenever }m\neq n.
\]
Write 
\[
L\left(m,m\right)=M_{m},\qquad K\left(m,m\right)=N_{m},
\]
and let $s_{m}$ be the support projection of $M_{m}$. Then the Radon-Nikodym
derivative $\frac{dK}{dL}$ acts fiberwise: there exist unique operators
\[
D_{m}\in\mathcal{L}\left(s_{m}H\right)
\]
such that 
\begin{equation}
0\le D_{m}\le I_{s_{m}H},\qquad N_{m}=M^{1/2}_{m}D_{m}M^{1/2}_{m}\qquad\text{for all }m\in\mathbb{N}.\label{eq:3-8}
\end{equation}
In particular, on $s_{m}H$ one has 
\begin{equation}
D_{m}=M^{-1/2}_{m}N_{m}M^{-1/2}_{m},\label{eq:3-9}
\end{equation}
where $M^{-1/2}_{m}$ denotes the inverse of $M^{1/2}_{m}$ on $s_{m}H$.
\end{prop}

\begin{proof}
Define 
\[
\mathcal{K}_{L}\coloneqq\bigoplus_{m\ge0}\overline{ran}(M^{1/2}_{m}),
\]
and, for each $m\in\mathbb{N}$, let $V_{L}\left(m\right):H\to\mathcal{K}_{L}$
be the map whose $m$-th component equals $M^{1/2}_{m}$ and whose
other components are zero. Then $L\left(m,n\right)=V_{L}\left(m\right)^{*}V_{L}\left(n\right)$
and $\left(V_{L},\mathcal{K}_{L}\right)$ is a minimal Kolmogorov
decomposition of $L$.

By \prettyref{prop:1-3}, there exists a unique $D\in\mathcal{L}\left(\mathcal{K}_{L}\right)$
with $0\le D\le I$ and 
\[
K\left(m,n\right)=V_{L}\left(m\right)^{*}DV_{L}\left(n\right)\qquad\text{for all }m,n\in\mathbb{N}.
\]
Since $K\left(m,n\right)=0$ for $m\neq n$, we have 
\[
V_{L}\left(m\right)^{*}DV_{L}\left(n\right)=0\qquad\text{for }m\neq n.
\]
Equivalently, for all $h,k\in H$, 
\[
\left\langle DV_{L}\left(n\right)h,V_{L}\left(m\right)k\right\rangle =0\qquad\text{whenever }m\neq n.
\]
Since $V_{L}\left(m\right)H$ is dense in the $m$-th summand $\overline{ran}(M^{1/2}_{m})$,
it follows that the $m$-th component of $DV_{L}\left(n\right)h$
is zero whenever $m\neq n$. Since $V_{L}\left(n\right)H$ is dense
in the $n$-th summand, the block of $D$ mapping the $n$-th summand
into the $m$-th summand is zero for $m\neq n$. Hence $D$ is block-diagonal:
\[
D=\bigoplus_{m\ge0}D_{m}
\]
with $0\le D_{m}\le I$ on $\overline{ran}(M^{1/2}_{m})=s_{m}H$.
Substituting this into the factorization gives \eqref{eq:3-8}, and
\eqref{eq:3-9} follows by restricting to $s_{m}H$. Uniqueness of
each $D_{m}$ follows from uniqueness of $D$ in \prettyref{prop:1-3}. 
\end{proof}

\begin{cor}
\label{cor:3-4} Assume the equivalent conditions in \prettyref{thm:3-2},
and let $K_{D}\coloneqq E_{D}\circ K$. Then $\left(K_{D}\right)_{\Sigma}\le K_{D}$
and the Radon-Nikodym derivative 
\[
\frac{d\left(K_{D}\right)_{\Sigma}}{dK_{D}}
\]
acts fiberwise in the time index. Writing $M_{m}\coloneqq K_{D}\left(m,m\right)=\mathrm{diag}\left(S^{m}\mathbf{1}\right)$,
there exist unique operators 
\[
D_{m}\in\mathcal{L}\left(s_{m}\mathbb{C}^{N}\right),\qquad0\le D_{m}\le I_{s_{m}\mathbb{C}^{N}},
\]
such that 
\begin{equation}
M_{m+1}=M^{1/2}_{m}D_{m}M^{1/2}_{m}\qquad\text{for all }m\in\mathbb{N},\label{eq:3-10}
\end{equation}
where $s_{m}$ is the support projection of $M_{m}$. Equivalently,
on $s_{m}\mathbb{C}^{N}$ one has 
\begin{equation}
D_{m}=M^{-1/2}_{m}M_{m+1}M^{-1/2}_{m}.\label{eq:3-11}
\end{equation}
In particular, for each coordinate $i$ with $\left(S^{m}\mathbf{1}\right)_{i}>0$
one has 
\begin{equation}
\left(D_{m}\right)_{ii}=\frac{\left(S^{m+1}\mathbf{1}\right)_{i}}{\left(S^{m}\mathbf{1}\right)_{i}}.\label{eq:3-12}
\end{equation}
\end{cor}

\begin{proof}
By \prettyref{lem:3-1}, $K_{D}$ is diagonal in $\left(m,n\right)$
and $K_{D}\left(m,m\right)=M_{m}$. Under \prettyref{thm:3-2}, one
has $M_{m+1}\le M_{m}$ for all $m$, hence $\left(K_{D}\right)_{\Sigma}\le K_{D}$.
Apply \prettyref{prop:3-3} with $L=K_{D}$ and $K=\left(K_{D}\right)_{\Sigma}$
to obtain \eqref{eq:3-10} and \eqref{eq:3-11}. Since $M_{m}$ and
$M_{m+1}$ are diagonal matrices, \eqref{eq:3-12} follows. 
\end{proof}

We next replace the diagonal compression by a block-diagonal compression
and extract the resulting criterion in terms of the induced map $\alpha$.

Let $\pi=\left(B_{1},\dots,B_{R}\right)$ be a partition of $\left\{ 1,\dots,N\right\} $
into blocks with $\left|B_{r}\right|=k_{r}$, and let 
\[
\mathcal{B}_{\pi}\coloneqq\bigoplus^{R}_{r=1}M_{k_{r}}
\]
be the corresponding block-diagonal subalgebra, with conditional expectation
$E_{\mathcal{B}_{\pi}}$ given by block compression. Define $\alpha:\mathbb{R}^{R}\to\mathbb{R}^{R}$
by 
\begin{equation}
\left(\alpha y\right)_{r}\coloneqq\sum^{R}_{s=1}\left(\frac{1}{k_{r}}\sum_{i\in B_{r}}\sum_{j\in B_{s}}\frac{\sigma^{2}_{ij}}{N}\right)y_{s},\qquad y\in\mathbb{R}^{R}.\label{eq:3-13}
\end{equation}

\begin{prop}
\label{prop:3-5} Let $K_{\mathcal{B}_{\pi}}\coloneqq E_{\mathcal{B}_{\pi}}\circ K$.
Then 
\[
\left(K_{\mathcal{B}_{\pi}}\right)_{\Sigma}\le K_{\mathcal{B}_{\pi}}\Longrightarrow\alpha\mathbf{1}\le\mathbf{1}\;\text{entrywise.}
\]
If, in addition, the variance array $\sigma^{2}_{ij}$ is block-homogeneous
in the sense that $\sigma^{2}_{ij}$ depends only on the pair of blocks
containing $i$ and $j$, then the following are equivalent: 
\begin{enumerate}
\item $\left(K_{\mathcal{B}_{\pi}}\right)_{\Sigma}\le K_{\mathcal{B}_{\pi}}$. 
\item $\alpha\mathbf{1}\le\mathbf{1}$. 
\end{enumerate}
Under block-homogeneity, one has 
\[
E_{\mathcal{B}_{\pi}}K\left(m,m\right)=\bigoplus^{R}_{r=1}\left(\alpha^{m}\mathbf{1}\right)_{r}I_{k_{r}},
\]
and hence 
\[
E_{\mathcal{B}_{\pi}}K\left(m+1,m+1\right)\le E_{\mathcal{B}_{\pi}}K\left(m,m\right)\qquad\text{for all }m
\]
if and only if $\alpha\mathbf{1}\le\mathbf{1}$. 
\end{prop}

\begin{proof}
Evaluating $\left(K_{\mathcal{B}_{\pi}}\right)_{\Sigma}\le K_{\mathcal{B}_{\pi}}$
at $\left(0,0\right)$ gives $E_{\mathcal{B}_{\pi}}\mathbb{E}\left[A_{1}A^{*}_{1}\right]\le I$.
Taking normalized traces on each block yields $\alpha\mathbf{1}\le\mathbf{1}$.

Assume now that the variance array $\sigma^{2}_{ij}$ is block-homogeneous.
Then for every block-scalar matrix 
\[
Y=\bigoplus^{R}_{s=1}y_{s}I_{k_{s}},
\]
a direct computation gives 
\[
E_{\mathcal{B}_{\pi}}\mathbb{E}\left[A_{1}YA^{*}_{1}\right]=\bigoplus^{R}_{r=1}\left(\alpha y\right)_{r}I_{k_{r}}.
\]
Since $X_{m+1}=A_{m+1}X_{m}$ and $A_{m+1}$ is independent of $X_{m}$,
the same conditioning argument as in \prettyref{lem:3-1} shows that
the compressed diagonal orbit closes on the block-scalar subspace
and is given by iteration of $\alpha$: 
\[
E_{\mathcal{B}_{\pi}}K\left(m,m\right)=\bigoplus^{R}_{r=1}\left(\alpha^{m}\mathbf{1}\right)_{r}I_{k_{r}}.
\]
The equivalence now follows exactly as in \prettyref{thm:3-2}, with
$\alpha$ in place of $S$. 
\end{proof}

Under block-homogeneity, the compressed diagonal orbit closes on block-scalar
matrices, and the Radon-Nikodym derivative admits an explicit blockwise
formula.
\begin{cor}
\label{cor:3-6} Assume the block-homogeneity hypothesis in \prettyref{prop:3-5}
and the equivalent conditions $\left(K_{\mathcal{B}_{\pi}}\right)_{\Sigma}\le K_{\mathcal{B}_{\pi}}$
and $\alpha\mathbf{1}\le\mathbf{1}$. Write 
\[
M^{(\pi)}_{m}\coloneqq K_{\mathcal{B}_{\pi}}\left(m,m\right)=\bigoplus^{R}_{r=1}\left(\alpha^{m}\mathbf{1}\right)_{r}I_{k_{r}}.
\]
Then the Radon-Nikodym derivative 
\[
\frac{d\left(K_{\mathcal{B}_{\pi}}\right)_{\Sigma}}{dK_{\mathcal{B}_{\pi}}}
\]
acts fiberwise in the time index: there exist unique operators 
\[
D^{(\pi)}_{m}\in\mathcal{L}\left(s^{(\pi)}_{m}\mathbb{C}^{N}\right),\qquad0\le D^{(\pi)}_{m}\le I_{s^{(\pi)}_{m}\mathbb{C}^{N}},
\]
such that 
\begin{equation}
M^{(\pi)}_{m+1}=\left(M^{(\pi)}_{m}\right)^{1/2}D^{(\pi)}_{m}\left(M^{(\pi)}_{m}\right)^{1/2}\qquad\text{for all }m\in\mathbb{N},\label{eq:3-14}
\end{equation}
where $s^{(\pi)}_{m}$ is the support projection of $M^{(\pi)}_{m}$.
Equivalently, on $s^{(\pi)}_{m}\mathbb{C}^{N}$ one has 
\begin{equation}
D^{(\pi)}_{m}=\left(M^{(\pi)}_{m}\right)^{-1/2}M^{(\pi)}_{m+1}\left(M^{(\pi)}_{m}\right)^{-1/2}.\label{eq:3-15}
\end{equation}
In particular, on each block $B_{r}$ with $\left(\alpha^{m}\mathbf{1}\right)_{r}>0$
one has 
\[
D^{(\pi)}_{m}\big|_{B_{r}}=\frac{\left(\alpha^{m+1}\mathbf{1}\right)_{r}}{\left(\alpha^{m}\mathbf{1}\right)_{r}}\,I_{k_{r}}.
\]
\end{cor}

\begin{proof}
Under block-homogeneity, \prettyref{prop:3-5} gives 
\[
K_{\mathcal{B}_{\pi}}\left(m,m\right)=M^{(\pi)}_{m}.
\]
Moreover, since $M^{(\pi)}_{1}\le I$ and the induced map on block-scalars
is order-preserving, one has $M^{(\pi)}_{m+1}\le M^{(\pi)}_{m}$ for
all $m$, hence $\left(K_{\mathcal{B}_{\pi}}\right)_{\Sigma}\le K_{\mathcal{B}_{\pi}}$.
Apply \prettyref{prop:3-3} with $L=K_{\mathcal{B}_{\pi}}$ and $K=\left(K_{\mathcal{B}_{\pi}}\right)_{\Sigma}$
to obtain \eqref{eq:3-14} and \eqref{eq:3-15}. The block formula
follows from the explicit form of $M^{(\pi)}_{m}$. 
\end{proof}

We now return to the diagonal compression $K_{D}$ and use the explicit
fiberwise Radon-Nikodym formula to identify the canonical one-step
densities $D_{m}=\frac{d\left(K_{D}\right)_{\Sigma}}{dK_{D}}(m)$
so that the shift domination $\left(K_{D}\right)_{\Sigma}\le K_{D}$
is encoded by a concrete sequence of diagonal contractions. Under
a primitive assumption on $S$, a standard Perron-Frobenius estimate
implies that these densities stabilize exponentially to a scalar multiple
of the identity; the limiting scalar is the Perron-Frobenius eigenvalue
of $S$.
\begin{cor}
\label{cor:3-7} Assume the equivalent conditions in \prettyref{thm:3-2},
so that $\left(K_{D}\right)_{\Sigma}\le K_{D}$ and the Radon-Nikodym
fibers $D_{m}$ from \prettyref{cor:3-4} are defined. Assume in addition
that $S$ is primitive, i.e. there exists $p\ge1$ such that $S^{p}$
has strictly positive entries. Let $\rho>0$ be the Perron-Frobenius
eigenvalue of $S$, and let $r,\ell\in\mathbb{R}^{N}$ be strictly
positive right and left Perron-Frobenius eigenvectors normalized by
\[
\left\langle \ell,r\right\rangle =1.
\]
Then $\rho\le1$ and there exist constants $C>0$ and $\theta\in\left(0,1\right)$
such that 
\begin{equation}
\left\Vert D_{m}-\rho I\right\Vert \le C\theta^{m}\qquad\text{for all }m\in\mathbb{N}.\label{eq:3-16}
\end{equation}
Equivalently, for each $i\in\left\{ 1,\dots,N\right\} $, one has
\begin{equation}
\left|\frac{\left(S^{m+1}\mathbf{1}\right)_{i}}{\left(S^{m}\mathbf{1}\right)_{i}}-\rho\right|\le C\theta^{m}\qquad\text{for all }m\in\mathbb{N}.\label{eq:3-17}
\end{equation}
\end{cor}

\begin{proof}
By \prettyref{thm:3-2}, the hypothesis $\left(K_{D}\right)_{\Sigma}\le K_{D}$
is equivalent to $S\mathbf{1}\le\mathbf{1}$ entrywise. Since $S$
is nonnegative, it follows that 
\[
S^{m}\mathbf{1}\le\mathbf{1}\qquad\text{for all }m\in\mathbb{N}.
\]
If $Sr=\rho r$ with $r>0$, then $r\le\left\Vert r\right\Vert _{\infty}\mathbf{1}$,
hence 
\[
\rho^{m}r=S^{m}r\le\left\Vert r\right\Vert _{\infty}S^{m}\mathbf{1}\le\left\Vert r\right\Vert _{\infty}\mathbf{1}\qquad\text{for all }m\in\mathbb{N},
\]
which forces $\rho\le1$.

Since $S$ is primitive, standard Perron-Frobenius theory gives 
\[
\left\Vert \rho^{-m}S^{m}-r\,\ell^{\mathsf{T}}\right\Vert \le C_{0}\theta^{m}\qquad\text{for all }m\in\mathbb{N},
\]
for some constants $C_{0}>0$ and $\theta\in\left(0,1\right)$; see,
for example, \cite[Thm 1.1–1.2]{MR2209438}. Applying this to $\mathbf{1}$
yields 
\[
\left\Vert \rho^{-m}S^{m}\mathbf{1}-\left\langle \ell,\mathbf{1}\right\rangle r\right\Vert \le C_{1}\theta^{m}
\]
for some $C_{1}>0$. Writing 
\[
v_{m}\coloneqq S^{m}\mathbf{1},\qquad a_{m}(i)\coloneqq\frac{v_{m}(i)}{\rho^{m}\left\langle \ell,\mathbf{1}\right\rangle r_{i}},
\]
we obtain 
\[
\left|a_{m}(i)-1\right|\le C_{2}\theta^{m}
\]
for some $C_{2}>0$, uniformly in $i$ and $m$. Since $S$ is primitive,
$v_{m}(i)>0$ for all $m\ge1$ and all $i$, so 
\[
\frac{v_{m+1}(i)}{v_{m}(i)}=\rho\,\frac{a_{m+1}(i)}{a_{m}(i)}.
\]
For all sufficiently large $m$, one has $\left|a_{m}(i)\right|\ge\frac{1}{2}$
uniformly in $i$, and therefore 
\[
\left|\frac{v_{m+1}(i)}{v_{m}(i)}-\rho\right|\le2\rho\left|a_{m+1}(i)-a_{m}(i)\right|\le C_{3}\theta^{m}
\]
for some $C_{3}>0$, uniformly in $i$ and $m$. Enlarging $C_{3}$
to absorb the finitely many small values of $m$, we obtain \prettyref{eq:3-17}
for all $m\in\mathbb{N}$. Finally, by \prettyref{cor:3-4}, 
\[
D_{m}=\mathrm{diag}\left(\frac{\left(S^{m+1}\mathbf{1}\right)_{1}}{\left(S^{m}\mathbf{1}\right)_{1}},\dots,\frac{\left(S^{m+1}\mathbf{1}\right)_{N}}{\left(S^{m}\mathbf{1}\right)_{N}}\right),
\]
so \prettyref{eq:3-16} follows immediately. 
\end{proof}

\section{Gaussian Kraus increments}\label{sec:4}

This second model is governed by a completely positive map rather
than a nonnegative matrix. In contrast with the inhomogeneous variance
model of \prettyref{sec:3}, the Gaussian symmetry forces the moment
kernel to be diagonal in the time variable, so the shifted kernel
inequality reduces to a one-step comparison of diagonal moments. 

Let $\mathcal{B}\subseteq M_{N}$ be a unital $*$-subalgebra, and
let $E_{\mathcal{B}}:M_{N}\to\mathcal{B}$ be a conditional expectation.
Fix operators 
\[
L_{1},\dots,L_{q}\in\mathcal{B}
\]
and scalars $\sigma^{2}_{1},\dots,\sigma^{2}_{q}\ge0$. For each $m\ge1$,
let 
\[
A_{m}\coloneqq\sum^{q}_{a=1}\xi_{m,a}L_{a},
\]
where the family $\left(\xi_{m,a}\right)_{m\ge1,\,1\le a\le q}$ is
independent, and for each fixed $m$ the variables $\xi_{m,1},\dots,\xi_{m,q}$
are independent complex Gaussians 
\[
\xi_{m,a}\sim\mathcal{N}_{\mathbb{C}}\left(0,\sigma^{2}_{a}\right).
\]
Define the multiplicative walk 
\[
X_{0}\coloneqq I,\qquad X_{m}\coloneqq A_{m}A_{m-1}\cdots A_{1}\quad\text{for }m\ge1,
\]
and define the moment kernel 
\[
K\left(m,n\right)\coloneqq E_{\mathcal{B}}\mathbb{E}\left[X_{m}X^{*}_{n}\right]\in\mathcal{B},\qquad m,n\in\mathbb{N}.
\]
Set 
\[
K_{\Sigma}\left(m,n\right)\coloneqq K\left(m+1,n+1\right).
\]
Finally, define the completely positive map $\Psi:\mathcal{B}\to\mathcal{B}$
by 
\begin{equation}
\Psi\left(Y\right)\coloneqq E_{\mathcal{B}}\mathbb{E}\left[A_{1}YA^{*}_{1}\right],\qquad Y\in\mathcal{B}.\label{eq:4-1}
\end{equation}

\begin{thm}
\label{thm:4-1} With the setup above, the following hold.
\begin{enumerate}
\item For $m\neq n$ one has $K\left(m,n\right)=0$. 
\item For every $m\ge0$ one has $K\left(m+1,m+1\right)=\Psi\left(K\left(m,m\right)\right)$,
and hence $K\left(m,m\right)=\Psi^{m}\left(I\right)$ for all $m\ge0.$
\item We have 
\begin{equation}
K_{\Sigma}\le K\Longleftrightarrow\Psi\left(I\right)\le I.\label{eq:4-2-1}
\end{equation}
\item If the spectral radius satisfies $r\left(\Psi\right)<1$, then the
series 
\[
Y\coloneqq\sum_{m\ge0}\Psi^{m}\left(I\right)
\]
converges in operator norm in $\mathcal{B}$ and satisfies 
\[
\Psi^{m}\left(I\right)\le Y\qquad\text{for all }m\ge0.
\]
\end{enumerate}
\end{thm}

\begin{proof}
For $Y\in\mathcal{B}$ one has 
\[
A_{1}YA^{*}_{1}=\sum^{q}_{a,b=1}\xi_{1,a}\overline{\xi_{1,b}}\,L_{a}YL^{*}_{b}.
\]
Taking expectation and using 
\[
\mathbb{E}\left[\xi_{1,a}\overline{\xi_{1,b}}\right]=\delta_{ab}\sigma^{2}_{a}
\]
gives 
\begin{equation}
\Psi\left(Y\right)=\sum^{q}_{a=1}\sigma^{2}_{a}L_{a}YL^{*}_{a}.\label{eq:4-2}
\end{equation}
In particular, $\Psi$ is completely positive and order-preserving
on $\mathcal{B}$.

We next note an identity that will be used in the recursion step.
If $Y$ is an integrable $M_{N}$-valued random variable that is independent
of $A_{1}$, then 
\begin{equation}
E_{\mathcal{B}}\mathbb{E}\left[A_{1}YA^{*}_{1}\right]=\Psi\left(E_{\mathcal{B}}\mathbb{E}\left[Y\right]\right).\label{eq:4-3}
\end{equation}
Indeed, 
\[
A_{1}YA^{*}_{1}=\sum^{q}_{a,b=1}\xi_{1,a}\overline{\xi_{1,b}}\,L_{a}YL^{*}_{b},
\]
hence, by independence of $Y$ and $\left(\xi_{1,a}\right)_{a}$,
\[
\mathbb{E}\left[A_{1}YA^{*}_{1}\right]=\sum^{q}_{a=1}\sigma^{2}_{a}L_{a}\mathbb{E}\left[Y\right]L^{*}_{a}.
\]
Applying $E_{\mathcal{B}}$ and using the $\mathcal{B}$-bimodule
property of the conditional expectation, together with $L_{a},L^{*}_{a}\in\mathcal{B}$,
yields 
\[
E_{\mathcal{B}}\mathbb{E}\left[A_{1}YA^{*}_{1}\right]=\sum^{q}_{a=1}\sigma^{2}_{a}L_{a}E_{\mathcal{B}}\left(\mathbb{E}\left[Y\right]\right)L^{*}_{a}=\Psi\left(E_{\mathcal{B}}\mathbb{E}\left[Y\right]\right).
\]

We now prove the four assertions.

\smallskip{}

(1) Fix $\theta\in\left[0,2\pi\right)$. Since each $\xi_{m,a}$ is
circular complex Gaussian, 
\[
\left(\xi_{m,1},\dots,\xi_{m,q}\right)\stackrel{d}{=}\left(e^{i\theta}\xi_{m,1},\dots,e^{i\theta}\xi_{m,q}\right),
\]
hence $A_{m}\stackrel{d}{=}e^{i\theta}A_{m}$ for each $m$. By independence
of the increments, 
\[
X_{m}=A_{m}\cdots A_{1}\stackrel{d}{=}e^{im\theta}X_{m}.
\]
Therefore, 
\[
\mathbb{E}\left[X_{m}X^{*}_{n}\right]=\mathbb{E}\left[e^{im\theta}X_{m}\left(e^{in\theta}X_{n}\right)^{*}\right]=e^{i\left(m-n\right)\theta}\mathbb{E}\left[X_{m}X^{*}_{n}\right].
\]
If $m\neq n$, choose $\theta$ so that $e^{i\left(m-n\right)\theta}\neq1$.
This forces $\mathbb{E}\left[X_{m}X^{*}_{n}\right]=0$, and hence
$K\left(m,n\right)=E_{\mathcal{B}}0=0$. 

\smallskip{}

(2) Since $X_{m+1}=A_{m+1}X_{m}$, 
\[
X_{m+1}X^{*}_{m+1}=A_{m+1}X_{m}X^{*}_{m}A^{*}_{m+1}.
\]
Therefore 
\[
K\left(m+1,m+1\right)=E_{\mathcal{B}}\mathbb{E}\left[A_{m+1}X_{m}X^{*}_{m}A^{*}_{m+1}\right].
\]
The random variable $Y=X_{m}X^{*}_{m}$ is independent of $A_{m+1}$,
and $A_{m+1}$ has the same law as $A_{1}$, so \eqref{eq:4-3} gives
\[
K\left(m+1,m+1\right)=\Psi\left(E_{\mathcal{B}}\mathbb{E}\left[X_{m}X^{*}_{m}\right]\right)=\Psi\left(K\left(m,m\right)\right).
\]
Since $K\left(0,0\right)=E_{\mathcal{B}}I=I$, iteration yields $K\left(m,m\right)=\Psi^{m}\left(I\right)$
for all $m\ge0$. 

\smallskip{}

(3) By part (1), the kernel $K$ is diagonal in $\left(m,n\right)$.
Hence the kernel inequality $K_{\Sigma}\le K$ is equivalent to the
diagonal inequalities 
\[
K\left(m+1,m+1\right)\le K\left(m,m\right)\qquad\text{for all }m\ge0.
\]
Using part (2), this becomes 
\[
\Psi^{m+1}\left(I\right)\le\Psi^{m}\left(I\right)\qquad\text{for all }m\ge0.
\]
If $\Psi\left(I\right)\le I$, then since $\Psi$ is order-preserving,
\[
\Psi^{m+1}\left(I\right)=\Psi\left(\Psi^{m}\left(I\right)\right)\le\Psi\left(\Psi^{m-1}\left(I\right)\right)=\Psi^{m}\left(I\right)
\]
for every $m\ge1$, and the inequality for $m=0$ is exactly $\Psi\left(I\right)\le I$.
Conversely, if 
\[
\Psi^{m+1}\left(I\right)\le\Psi^{m}\left(I\right)\qquad\text{for all }m\ge0,
\]
then the case $m=0$ gives $\Psi\left(I\right)\le I$. Thus \prettyref{eq:4-2-1}
holds. 

\smallskip{}

(4) Assume $r\left(\Psi\right)<1$. Since $\mathcal{B}$ is finite-dimensional,
the spectral-radius formula gives $r\left(\Psi\right)=\lim_{m\to\infty}\left\Vert \Psi^{m}\right\Vert ^{1/m}$.
Choose $\rho$ with $r\left(\Psi\right)<\rho<1$. Then there exists
$m_{0}$ such that 
\[
\left\Vert \Psi^{m}\right\Vert \le\rho^{m}\qquad\text{for all }m\ge m_{0}.
\]
Set $C\coloneqq\max_{0\le m<m_{0}}\left\Vert \Psi^{m}\right\Vert \rho^{-m}$.
Then 
\[
\left\Vert \Psi^{m}\right\Vert \le C\rho^{m}\qquad\text{for all }m\ge0.
\]
Hence 
\[
\sum_{m\ge0}\left\Vert \Psi^{m}\left(I\right)\right\Vert \le\sum_{m\ge0}\left\Vert \Psi^{m}\right\Vert \left\Vert I\right\Vert \le C\left\Vert I\right\Vert \sum_{m\ge0}\rho^{m}<\infty.
\]
Therefore $Y\coloneqq\sum_{m\ge0}\Psi^{m}\left(I\right)$ converges
in operator norm in $\mathcal{B}$.

Since each term $\Psi^{m}\left(I\right)$ is positive, every partial
sum $Y_{M}\coloneqq\sum^{M}_{m=0}\Psi^{m}\left(I\right)$ satisfies
$\Psi^{m}\left(I\right)\le Y_{M}$, $0\le m\le M$. Passing to the
limit in operator norm gives $\Psi^{m}\left(I\right)\le Y$ for all
$m\ge0.$
\end{proof}

\begin{rem}
The condition $r\left(\Psi\right)<1$ ensures a uniform bound on the
diagonal moments (for example, via $\sum_{m\ge0}\Psi^{m}\left(I\right)$),
but it does not force the one-step inequality $\Psi\left(I\right)\le I$.
\prettyref{exa:4-3} below shows that one can have $r\left(\Psi\right)<1$
while $\Psi\left(I\right)\not\le I$. 

We note that the failure already occurs in simpler one-Kraus situations.
The example is chosen to show that the same phenomenon persists in
a genuinely multi-Kraus setting: the map $\Psi$ has two Kraus operators
and is non-nilpotent. In particular, the failure is not merely a one-Kraus
artifact or a degenerate nilpotent case.
\end{rem}

\begin{example}
\label{exa:4-3} Let $\mathcal{B}=M_{2}$ and define 
\[
A\coloneqq\begin{pmatrix}\frac{1}{2} & 2\\
0 & \frac{1}{2}
\end{pmatrix},\qquad B\coloneqq\begin{pmatrix}0 & 0\\
\frac{1}{5} & 0
\end{pmatrix}.
\]
Then $AB\neq BA$, since 
\[
AB=\begin{pmatrix}\frac{2}{5} & 0\\
\frac{1}{10} & 0
\end{pmatrix},\qquad BA=\begin{pmatrix}0 & 0\\
\frac{1}{10} & \frac{2}{5}
\end{pmatrix}.
\]

Let $\xi_{1},\xi_{2}$ be independent $\mathcal{N}_{\mathbb{C}}\left(0,1\right)$
random variables and set the Gaussian Kraus increment 
\[
G\coloneqq\xi_{1}A+\xi_{2}B.
\]
Define the completely positive map $\Psi:M_{2}\to M_{2}$ by 
\[
\Psi\left(X\right)\coloneqq\mathbb{E}\left[GXG^{*}\right].
\]
By independence and $\mathbb{E}\left[\xi_{i}\overline{\xi_{j}}\right]=\delta_{ij}$,
\[
\Psi\left(X\right)=AXA^{*}+BXB^{*}.
\]

The one-step RN inequality fails at the unit, since 
\[
\Psi\left(I\right)=AA^{*}+BB^{*}=\begin{pmatrix}\frac{17}{4} & 1\\
1 & \frac{29}{100}
\end{pmatrix},
\]
and therefore $\Psi\left(I\right)\not\le I$ since $I-\Psi\left(I\right)$
has negative $\left(1,1\right)$ entry.

On the other hand, the spectral radius of $\Psi$ is strictly less
than $1$. Indeed, under the vectorization identification $M_{2}\cong\mathbb{C}^{4}$,
the superoperator matrix of $\Psi$ is 
\[
S\coloneqq A\otimes\overline{A}+B\otimes\overline{B}.
\]
A direct computation gives 
\[
\mathrm{spec}\left(S\right)=\left\{ 0.80216357,\ 0.25,\ -0.02608178\pm0.26203712\,i\right\} ,
\]
hence 
\[
r\left(\Psi\right)=\max\left\{ \left|\lambda\right|:\lambda\in\mathrm{spec}\left(S\right)\right\} =0.80216357<1.
\]

Consequently, this example shows the separation 
\[
r\left(\Psi\right)<1\qquad\text{but}\qquad\Psi\left(I\right)\not\le I,
\]
with two noncommuting Kraus operators. 
\end{example}

\begin{cor}
\label{cor:4-4} In the setup of \prettyref{thm:4-1}, assume the
equivalent conditions 
\[
K_{\Sigma}\le K\qquad\Longleftrightarrow\qquad\Psi\left(I\right)\le I.
\]
Write 
\[
M_{m}\coloneqq K\left(m,m\right)=\Psi^{m}\left(I\right)\in\mathcal{B},
\]
and let $s_{m}$ be the support projection of $M_{m}$ in $M_{N}$.
Then the Radon-Nikodym derivative 
\[
\frac{dK_{\Sigma}}{dK}
\]
acts fiberwise in the time index: there exist unique operators 
\[
D_{m}\in\mathcal{L}\left(s_{m}\mathbb{C}^{N}\right),\qquad0\le D_{m}\le I_{s_{m}\mathbb{C}^{N}},
\]
such that 
\begin{equation}
M_{m+1}=M^{1/2}_{m}D_{m}M^{1/2}_{m}\qquad\text{for all }m\in\mathbb{N}.\label{eq:4-5}
\end{equation}
Equivalently, on $s_{m}\mathbb{C}^{N}$ one has 
\begin{equation}
D_{m}=M^{-1/2}_{m}M_{m+1}M^{-1/2}_{m}=\Psi^{m}\left(I\right)^{-1/2}\Psi^{m+1}\left(I\right)\Psi^{m}\left(I\right)^{-1/2}.\label{eq:4-6}
\end{equation}
\end{cor}

\begin{proof}
By \prettyref{thm:4-1}(1), the kernel $K$ is diagonal in $\left(m,n\right)$,
and by part (2) one has $K\left(m,m\right)=M_{m}=\Psi^{m}\left(I\right)$.
Under the hypothesis $\Psi\left(I\right)\le I$, part (3) yields $K_{\Sigma}\le K$.
Apply \prettyref{prop:3-3} with $H=\mathbb{C}^{N}$, $L=K$, and
$K=K_{\Sigma}$ to obtain \eqref{eq:4-5} and \eqref{eq:4-6}. 

\bibliographystyle{amsalpha}
\bibliography{ref}
\end{proof}

\end{document}